\documentclass[11pt]{article}

\newtheorem{theo}{Theorem}

\newtheorem{lemma}{Lemma}[section]
\newtheorem{definition}[lemma]{Definition}

\newcommand{\ignore}[1]{}
\def\square{\vrule height6pt width7pt depth1pt}
\def\endpf{\hfill\square\bigskip}

\begin{document}
\title{Quasi-randomness is determined by the distribution of copies of a fixed graph in 
equicardinal large sets}
\author{
\and
Raphael Yuster
\thanks{Department of Mathematics, University of Haifa, Haifa
31905, Israel. E--mail: raphy@math.haifa.ac.il} }

\date{}
\maketitle

\setcounter{page}{1}
\begin{abstract}

For every fixed graph $H$ and every fixed $0 < \alpha < 1$, we show that if a graph $G$ has 
the property that all subsets of size $\alpha n$
contain the ``correct'' number of copies of $H$ one would expect to
find in the random graph $G(n,p)$ then $G$ behaves like the random graph
$G(n,p)$; that is, it is $p$-quasi-random in the sense of Chung, Graham, and Wilson 
\cite{ChGrWi}. This solves a conjecture raised by Shapira \cite{Sh} and solves in a strong 
sense an open problem of Simonovits and S\'os \cite{SiSo}.

\end{abstract}

\thispagestyle{empty}
\setcounter{page}{1}

\section{Introduction}
The theory of quasi-random graphs asks the following fundamental question: which properties 
of graphs are such that any graph that satisfies them, resembles an appropriate random graph 
(namely, the graph satisfies the properties that a random graph would
satisfy, with high probability). Such properties are called {\em quasi-random}.

The theory of quasi-random graphs was initiated by Thomason \cite{Th1,Th2} and then
followed by Chung, Graham and Wilson who proved the fundamental theorem of quasi-random 
graphs \cite{ChGrWi}. Since then there have been many papers on this subject (see, e.g.
the excellent survey \cite{KrSu}).
Quasi-random properties were also studied for
other combinatorial structures such as set systems \cite{ChGr}, tournaments \cite{ChGr1}, and 
hypergraphs \cite{ChGr2}. There are also some very recent results on quasi-random groups
\cite{Go} and generalized quasi-random graphs \cite{LoSo}.

In order to formally define $p$-quasi-randomness we need to state the fundamental theorem of 
quasi-random graphs. As usual, a {\em labeled copy} of a graph $H$ in a graph 
$G$ is an injective mapping $\phi$ from $V(H)$ to $V(G)$ that maps edges to edges.
That is $(x,y) \in E(H)$ implies $(\phi(x), \phi(y)) \in E(G)$.
For a set of vertices $U \subset V(G)$ we denote by $H[U]$ the number of labeled
copies of $H$ in the subgraph of $G$ induced by $U$ and by $e(U)$ the number of edges
of $G$ with both endpoints in $U$. A graph sequence $(G_n)$ is an infinite sequence of graphs
$\{G_1,G_2,\ldots\}$ where $G_n$ has $n$ vertices. The following
result of Chung, Graham, and Wilson \cite{ChGrWi}
shows that many properties of different nature are equivalent to the notion of 
quasi-randomness, defined using edge distribution.
The original theorem lists seven such equivalent properties, but we only state four of them 
here.

\begin{theo}[Chung, Graham, and Wilson \cite{ChGrWi}]
\label{t-CGW}
Fix any $1 < p < 1$.
For any graph sequence $(G_n)$ the following properties are equivalent:
\begin{itemize}
\item[${\cal P}_1(t)$:]
For an even integer $t \ge 4$, let
$C_t$ denote the cycle of length $t$. Then $e(G_n) = \frac{1}{2}pn^2+o(n^2)$ and
$C_t[G_n] = p^tn^t + o(n^t)$.
\item[${\cal P}_2$:]
For any subset of vertices $U \subseteq V(G_n)$ we have
$e(U)=\frac{1}{2}p|U|^2 +o(n^2)$.
\item[${\cal P}_3$:]
For any subset of vertices $U \subseteq V(G_n)$ of size $n/2$ we have
$e(U)=\frac{1}{2}p|U|^2 +o(n^2)$.
\item[${\cal P}_4(\alpha)$:]
Fix an $\alpha \in (0,\frac12)$. For any
$U \subseteq V(G_n)$ of size $\alpha n$ we have $e(U, V \setminus U)=p\alpha(1-\alpha)n^2 
+o(n^2)$.
\end{itemize}
\end{theo}
The {\em formal} meaning of the properties being equivalent is expressed, as
usual, using $\epsilon,\delta$ notation. For example the meaning that ${\cal P}_3$ implies
${\cal P}_2$ is that for any $\epsilon > 0$ there exist $\delta=\delta(\epsilon)$
and $N=N(\epsilon)$ so that for all $n > N$, if $G$ is a graph with $n$ vertices
having the property that any subset of vertices $U$ of size $n/2$ satisfies
$|e(U)-\frac{1}{2}p|U|^2| < \delta n^2$ then also for any subset of vertices $W$
we have $|e(W)-\frac{1}{2}p|W|^2| < \epsilon n^2$.

Given Theorem \ref{t-CGW} we say that a graph property is $p$-quasi-random if it is 
equivalent to any (and therefore all) of the four properties defined in that theorem.
(We will usually just say {\em quasi-random} instead of {\em $p$-quasi-random} since $p$
is fixed throughout the proofs).
Note, that each of the four properties in Theorem \ref{t-CGW} is a property we would expect 
$G(n,p)$ to satisfy with high probability.

It is far from true, however, that any property that almost surely holds for $G(n,p)$ is 
quasi-random. For example, it is easy to see that having vertex degrees
$np(1+o(1))$ is not a quasi-random property (just take vertex-disjoint cliques of 
size roughly $np$ each). An important family of {\em non} quasi-random properties are those requiring the 
graphs in the sequence to have the correct number of copies of a fixed graph $H$. Note that 
${\cal P}_1(t)$ guarantees that for any {\em even} $t$, if a graph sequence
has the correct number of edges as well as the correct number of copies of $H=C_t$ then the 
sequence is quasi-random. As observed in \cite{ChGrWi} this is not true for all graphs.
In fact, already for $H=K_3$ there are simple constructions showing that this is not true.

Simonovits and S\'os observed that the standard counter-examples showing that for some graphs 
$H$, having the correct number of copies of $H$ is not enough to guarantee quasi-randomness, 
have the property that the number of copies of $H$ in some of the induced subgraphs of these 
counter-examples deviates significantly from what it should be.
As quasi-randomness is a hereditary property, in the sense that we expect a
sub-structure of a random-like object to be random-like as well, they introduced the 
following variant
of property ${\cal P}_1$ of Theorem \ref{t-CGW}, where now we require all subsets of vertices 
to contains the ``correct'' number of copies of $H$.

\begin{definition}[${\cal P}_H$]
\label{d-ss}
For a fixed graph $H$ with $h$ vertices and $r$ edges,
we say that a graph sequence $(G_n)$ satisfies ${\cal P}_H$ if
all subsets of vertices $U \subset V(G_n)$ with $|U|=\alpha n$
satisfy $H[U] = p^r|U|^h+o(n^h)$.
\end{definition}

As opposed to ${\cal P}_1$, which is quasi-random only for even cycles, Simonovits
and S\'os \cite{SiSo} showed that ${\cal P}_H$ is quasi-random for any nonempty graph $H$.
\begin{theo}
\label{t-ss}
For any fixed $H$ that has edges, property ${\cal P}_H$ is quasi-random.
\end{theo}

We can view property ${\cal P}_H$ as a generalization of property ${\cal P}_2$ in
Theorem \ref{t-CGW}, since ${\cal P}_2$ is just the special case ${\cal P}_{K_2}$.
Now, property ${\cal P}_3$ in Theorem \ref{t-CGW} guarantees that in order
to infer that a sequence is quasi-random, and thus satisfies ${\cal P}_2$,
it is enough to require only the sets
of vertices of size $n/2$ to contain the correct number of edges.
An open problem raised by Simonovits and S\'os \cite{SiSo}, and in a stronger form by Shapira 
\cite{Sh}, is that the analogous condition also holds for any $H$. Namely, in order
to infer that a sequence is quasi-random, and thus satisfies ${\cal P}_H$,
it is enough, say, to require only the sets of vertices of size $n/2$ to contain the correct 
number of copies of $H$. Shapira \cite{Sh} proved that is it enough to consider sets of 
vertices of size
$n/(h+1)$. Hence, in his result, the cardinality of the sets {\em depends} on $h$.
Thus, if $H$ has 1000 vertices, Shpiras's result shows that it suffices to check vertex
subsets having a fraction smaller than $1/1000$ of the total number of vertices.
His proof method cannot be extended to obtain the same result for fractions larger than 
$1/(h+\epsilon)$.

In this paper we settle the above mentioned open problem completely. In fact, we show that 
for any $H$,
not only is it enough to check only subsets of size $n/2$, but, more generally, we show that
it is enough to check subsets of size $\alpha n$ for any fixed $\alpha \in (0,1)$.
More formally, we define:

\begin{definition}[${\cal P}_{H,\alpha}$]
\label{d-main}
For a fixed graph $H$ with $h$ vertices and $r$ edges and fixed $0 < \alpha < 1$
we say that a graph sequence $(G_n)$ satisfies ${\cal P}_{H,\alpha}$ if
all subsets of vertices $U \subset V(G_n)$ with $|U|=\lfloor \alpha n \rfloor$
satisfy $H[U] = p^r|U|^h+o(n^h)$.
\end{definition}

\noindent
Our main result is, therefore:

\begin{theo}
\label{t-main}
For any fixed graph $H$ and any fixed $0 < \alpha < 1$, property ${\cal P}_{H, \alpha}$ is 
quasi-random.
\end{theo}

\section{Proof of the main result}

For the remainder of this section let $H$ be a fixed graph with $h > 1$ vertices
and $r > 0$ edges, and let $\alpha \in (0,1)$ be fixed.
Throughout this section we ignore rounding issues and, in particular, assume that $\alpha n$
is an integer, as this has no effect on the asymptotic nature of the results.
 
Suppose that the graph sequence $(G_n)$ satisfies ${\cal P}_{H,\alpha}$.
We will prove that it is quasi-random by showing that it also satisfies
${\cal P}_H$. In other words, we need to prove the following lemma which,
together with Theorem \ref{t-ss}, yields Theorem \ref{t-main}.
\begin{lemma}
\label{l-main}
For any $\epsilon > 0$ there exists $N=N(\epsilon,h,\alpha)$ and 
$\delta=\delta(\epsilon,h,\alpha)$ so that for all $n > N$,
if $G$ is a graph with $n$ vertices satisfying that for all $U \subset V(G)$
with $|U|=\alpha n$ we have $|H[U] - p^r|U|^h| < \delta n^h$
then $G$ also satisfies that for all $W \subset V(G)$
we have $|H[W] - p^r|W|^h| < \epsilon n^h$.
\end{lemma}
{\bf Proof:}\,
Suppose therefore that $\epsilon > 0$ is given. Let
$N=N(\epsilon,h,\alpha)$, $\epsilon'=\epsilon'(\epsilon,h,\alpha)$ and
$\delta=\delta(\epsilon,h,\alpha)$ be parameters
to be chosen so that $N$ is sufficiently large and $\delta \ll \epsilon'$ are both 
sufficiently
small to satisfy the inequalities that will follow, and it will be clear that
they are indeed only functions of $\epsilon,h$, and $\alpha$.

Now, let $G$ be a graph with $n > N$ vertices satisfying that
for all $U \subset V(G)$ with $|U|=\alpha n$ we have $|H[U] - p^r|U|^h| < \delta n^h$.
Consider any subset $W \subset V(G)$. We need to prove that
$|H[W] - p^r|W|^h| < \epsilon n^h$.

For convenience, set $k=\alpha n$.
Let us first prove this for the case where $|W|=m > k$.
This case can rather easily be proved via a simple counting argument.
Denote by ${\cal U}$ the set of ${m \choose k}$ $k$-subsets of $W$.
Hence, by the given condition on $k$-subsets,
\begin{equation}
\label{e1}
{m \choose k}(p^rk^h - \delta n^h) < \sum_{U \in {\cal U}} H[U] < {m \choose k}(p^rk^h + 
\delta n^h) \,.
\end{equation}
Every copy of $H$ in $W$ appears in precisely ${{m-h} \choose {k-h}}$ distinct
$U \in {\cal U}$. it follows from (\ref{e1}) that
\begin{equation}
\label{e2}
H[W] = \frac{1}{{{m-h} \choose {k-h}}}\sum_{U \in {\cal U}} H[U] <
\frac{{m \choose k}}{{{m-h} \choose {k-h}}}(p^rk^h + \delta n^h) <
p^rm^h + \frac{\epsilon'}{2} n^h \,,
\end{equation}
and similarly from (\ref{e1})
\begin{equation}
\label{e3}
H[W] = \frac{1}{{{m-h} \choose {k-h}}}\sum_{U \in {\cal U}} H[U] >
\frac{{m \choose k}}{{{m-h} \choose {k-h}}}(p^rk^h - \delta n^h) >
p^rm^h - \frac{\epsilon'}{2} n^h \,.
\end{equation}

We now consider the case where $|W| = m = \beta n < \alpha n = k$.
Notice that we can assume that $\beta \ge \epsilon$ since otherwise the
result is trivially true.
The set ${\cal H}$ of $H$-subgraphs of $G$ can be partitioned into
$h+1$ types, according to the number of vertices they have in $W$.
Hence, for $j=0,\ldots,h$ let ${\cal H}_j$ be the set of $H$-subgraphs of
$G$ that contain precisely $j$ vertices in $V \setminus W$. Notice that, by definition,
$|{\cal H}_0| = H[W]$. For convenience, denote $w_j = |{\cal H}_j|/n^h$.
We therefore have
\begin{equation}
\label{e4}
w_0 + w_1 + \cdots + w_h = \frac{|{\cal H}|}{n^h} = \frac{H[V]}{n^h} = p^r+\mu
\end{equation}
where $|\mu| < \epsilon'/2$.

Define $\lambda=\frac{(1-\alpha)}{h+1}$ and set $k_i=k+i\lambda n$ for
$i=1,\ldots,h$. Let $Y_i \subset V \setminus W$ be a random set of $k_i-m$
vertices, chosen uniformly at random from all ${{n-m} \choose {k_i-m}}$ subsets
of size $k_i-m$ of $V \setminus W$. Denote $K_i=Y_i \cup W$ and notice that $|K_i| = k_i > 
\alpha n$.
We will now estimate the number of elements of ${\cal H}_j$ that ``survive'' in $K_i$.
Formally, let ${\cal H}_{j,i}$ be the set of elements of ${\cal H}_j$ that have all of their
vertices in $K_i$, and let $m_{j,i}=|{\cal H}_{j,i}|$. Clearly, $m_{0,i}=H[W]$ since
$W \subset K_i$. Furthermore, by (\ref{e2}) and (\ref{e3}),
\begin{equation}
\label{e5}
m_{0,i}+m_{1,i}+\cdots+m_{h,i} = H[K_i] = p^rk_i^h+\rho_in^h
\end{equation}
where $\rho_i$ is a random variable with $|\rho_i| < \epsilon'/2$.

For an $H$-copy
$T \in {\cal H}_j$ we compute the probability $p_{j,i}$ that $T \in H[K_i]$.
Since $T \in H[K_i]$ if and only if all the $j$ vertices of $T$ in $V \setminus W$
appear in $Y_i$ we have
$$
p_{j,i} = \frac{{{n-m-j} \choose {k_i-m-j}}}{{{n-m} \choose {k_i-m}}} =
\frac{(k_i-m)\cdots(k_i-m-j+1)}{(n-m) \cdots (n-m-j+1)}\,.
$$
Defining $x_i = (k_i-m)/(n-m)$ and noticing that
$$
x_i = \frac{k_i-m}{n-m} = \frac{\alpha - \beta}{1-\beta}+\frac{\lambda}{1-\beta}i
$$
it follows that
\begin{equation}
\label{e6}
\left|p_{j,i}-x_i^j\right| < \frac{\epsilon'}{2}\,.
\end{equation}

Clearly, the expectation of $m_{j,i}$ is ${\rm E}[m_{j,i}]=p_{j,i}|{\cal H}_j|$.
By linearity of expectation we have from (\ref{e5}) that
$$
{\rm E}[m_{0,i}]+{\rm E}[m_{1,i}]+\cdots+{\rm E}[m_{h,i}] =
{\rm E}[H[K_i]] = p^rk_i^h+{\rm E}[\rho_i]n^h.
$$
Dividing the last equality by $n^h$ we obtain
\begin{equation}
\label{e7}
p_{0,i}w_0 + \cdots + p_{h,i}w_h = p^r\left(\alpha+\lambda i\right)^h+E[\rho_i]\,.
\end{equation}
By (\ref{e6}) and (\ref{e7}) we therefore have
\begin{equation}
\label{e8}
\sum_{j=0}^h x_i^j w_j = p^r\left(\alpha+\lambda i\right)^h+\mu_i
\end{equation}
where $\mu_i = E[\rho_i]+\zeta_i$ and $|\zeta_i| < \epsilon'/2$.
Since also $|\rho_i| < \epsilon'/2$ we have that $|\mu_i| < \epsilon'$.

Now, (\ref{e4}) and (\ref{e8}) form together a system of $h+1$ linear equations with
the $h+1$ variables $w_0,\ldots,w_h$. The coefficient matrix of this system
is just the Vandermonde matrix $A=A(x_1,\ldots,x_h,1)$.
Since $x_1,\ldots,x_h,1$ are all distinct, and, in fact, the gap between any two of them
is at least $\lambda/(1-\beta)=(1-\alpha)/((h+1)(1-\beta)) \ge (1-\alpha)/(h+1)$,
we have that the system has a unique solution which is $A^{-1}b$ where
$b \in R^{h+1}$ is the column vector whose $i$'th coordinate is $p^r\left(\alpha+\lambda 
i\right)^h+\mu_i$
for $i=1,\ldots,h$ and whose last coordinate is $p^r + \mu$.
Consider now the vector $b^*$ which is the same as $b$, just without the $\mu_i$'s. Namely
$b^* \in R^{h+1}$ is the column vector whose $i$'th coordinate is $p^r\left(\alpha+\lambda 
i\right)^h$
for $i=1,\ldots,h$ and whose last coordinate is $p^r$.
Then the system $A^{-1}b^*$ also has a unique solution and, in fact, we {\em know} explicitly 
what this solution
is. It is the vector $w^*=(w_0^*,\ldots,w^*_h)$ where
$$
w_j^* = p^r {h \choose j}\beta^{h-j}(1-\beta)^j\,.
$$
Indeed, it is straightforward to verify the equality
$$
\sum_{j=0}^h 
p^r {h \choose j}\beta^{h-j}(1-\beta)^j = p^r
$$
and, for all $i=1,\ldots,h$ the equalities
$$
\sum_{j=0}^h \left(\frac{\alpha - \beta}{1-\beta}+\frac{\lambda}{1-\beta}i\right)^j
p^r {h \choose j}\beta^{h-j}(1-\beta)^j = p^r\left(\alpha+\lambda i\right)^h\,.
$$
Now, since the mapping $F: R^{h+1} \rightarrow R^{h+1}$ mapping a vector $c$ to $A^{-1}c$ is 
continuous,
we know that for $\epsilon'$ sufficiently small, if each coordinate of $c$ has absolute value
less than $\epsilon'$, then each coordinate of $A^{-1}c$ has absolute value at most 
$\epsilon$.
Now, define $c=b-b^*=(\mu_1,\ldots,\mu_h,\mu)$. Then we have that each coordinate $w_i$ of 
$A^{-1}b$
differs from the corresponding coordinate $w_i^*$ of $A^{-1}b^*$ by at most $\epsilon$.
In particular,
$$
|w_0 -w_0^*| = |w_0 - p^r\beta^h| < \epsilon.
$$
Hence,
$$
|H[W] - n^hp^r\beta^h | = |H[W] - p^r|W|^h| < \epsilon n^h
$$
as required.
\endpf

\end{document}